\documentclass{amsart}

\newtheorem{theorem}{Theorem}[section]
\newtheorem{lemma}[theorem]{Lemma}
\newtheorem{proposition}[theorem]{Proposition}
\newtheorem{corollary}[theorem]{Corollary}

\theoremstyle{definition}

\newtheorem{remark}[theorem]{Remark}
\newtheorem{question}[theorem]{Question}
\newtheorem*{acknowledgements}{Acknowledgements}

\theoremstyle{remark}

\usepackage{amscd,amssymb}

\newcommand\mylabel[1]{\label{#1}}

\newcommand{\ZZ}{\mathbb{Z}}
\newcommand{\QQ}{\mathbb{Q}}

\newcommand{\PP}{\mathbb{P}}
\renewcommand{\AA}{\mathbb{A}}
\newcommand{\GG}{\mathbb{G}}

\newcommand  {\shA}     {\mathcal{A}}
\newcommand  {\shAut}   {\mathcal{A}\!\text{\textit{ut}}}
\newcommand  {\shB}     {\mathcal{B}}

\newcommand  {\shDiv}   {\mathcal{D} \!\text{\textit{iv}}}

\newcommand  {\shEnd}   {\mathcal{E}\!\text{\textit{nd}}}

\newcommand  {\shE}     {\mathcal{E}}
\newcommand  {\shF}     {\mathcal{F}}
\newcommand  {\shG}     {\mathcal{G}}
\newcommand  {\shH}     {\mathcal{H}}

\newcommand  {\shI}     {\mathcal{I}}

\newcommand  {\shM}     {\mathcal{M}}

\newcommand  {\shL}     {\mathcal{L}}

\newcommand  {\shS}     {\mathcal{S}}
\newcommand  {\shT}     {\mathcal{T}}

\newcommand  {\shP}     {\mathcal{P}}

\newcommand  {\shZ}     {\mathcal{Z}}

\newcommand  {\bBr}     {\widetilde{\operatorname{Br}}}
\newcommand  {\Br}      {\operatorname{Br}}

\newcommand  {\Cl}      {\operatorname{Cl}}
\newcommand  {\cH}      {\check{H}}

\newcommand  {\Div}     {\operatorname{Div}}

\newcommand  {\et}      {{\text{\rm \'{e}t}}}

\newcommand  {\Gr}      {\operatorname{Gr}}

\newcommand  {\id}      {\operatorname{id}}

\newcommand  {\Isom}    {\operatorname{Isom}}

\newcommand  {\dirlim}  {\varinjlim}
\newcommand  {\invlim}  {\varprojlim}

\newcommand  {\lra}     {\longrightarrow}

\renewcommand{\O}       {\mathcal{O}}

\newcommand  {\Pic}     {\operatorname{Pic}}

\newcommand  {\quadand} {\quad\text{and}\quad}

\newcommand  {\ra}      {\rightarrow}

\newcommand  {\Spec}    {\operatorname{Spec}}
\newcommand  {\sh}      {{\operatorname{sh}}}

\newcommand {\zar}      {{\operatorname{zar}}}

\def\mydate{\number\day\space\ifcase\month \or January\or February\or March\or 
April\or May\or June\or July\or
August\or September\or October\or November\or December\fi \space\number\year}

\begin{document}

\title[The bigger Brauer group]{The bigger Brauer group is really big}

\author[Stefan Schroer]{Stefan Schr\"oer}
\address{Mathematische Fakult\"at, Ruhr-Universit\"at, 
         44780 Bochum, Germany}
\email{s.schroeer@ruhr-uni-bochum.de}


\dedicatory{final version, 19 July 2002}

\begin{abstract}
I show that each \'etale $n$-cohomology class on  noetherian
schemes comes from a \v{C}ech cocycle, provided that
any $n$-tuple of points  admits an affine open neighborhood.
Together with results of 
Raeburn and Taylor on the bigger Brauer group,
this implies  that for schemes such 
that each pair of points admits an affine open neighborhood,
any \'etale $\GG_m$-gerbe comes
from a coherent central separable algebra.
Such algebras are nonunital generalizations of Azumaya algebras.
I also prove that, on normal noetherian schemes,
each Zariski $\GG_m$-gerbe comes from a central separable algebra.
\end{abstract}

\maketitle

\section*{Introduction}

Grothendieck \cite{GB} asked whether each torsion class 
in $H^2_\et(X,\GG_m)$ on a  scheme $X$
comes from an Azumaya algebra.
This is a major open 
problem in the theory of Brauer groups.
Gabber \cite{Gabber 1981} proved it  for affine schemes.
But even for smooth projective threefolds the answer seems
to be unknown. Edidin, Hassett, Kresch, and Vistoli
\cite{Edidin; Hassett; Kresch; Vistoli 1999} recently
found a counterexamples for nonseparated schemes.

To attack the problem, it is perhaps a good idea to modify it.
Taylor \cite{Taylor 1982}  generalized the notion of Azumaya algebras to
\emph{central separable
algebras}, which 
are not necessarily locally free or unital.
Nevertheless, they come along with a  $\GG_m$-gerbe of splittings 
and therefore define a cohomology class in $H^2_\et(X,\GG_m)$.
Assuming that each finite subset in $X$ admits an affine open
neighborhood, Raeburn and Taylor \cite{Raeburn; Taylor 1985}
proved that each 2-cohomology class, torsion or not,
comes from a coherent central separable algebra.
Caenepeel and Grandjean \cite{Caenepeel; Grandjean 1998} 
later fixed some problems in 
the original arguments.

Actually, the arguments of Raeburn and Taylor show that,
on arbitrary noetherian schemes, 
each \v{C}ech 2-cohomology class comes from a coherent
central separable algebra. 
Not every 2-cohomology class, however, comes from \v{C}ech cocycles.
Rather, the obstruction  is
a 1-cocycle class with values in the presheaf $U\mapsto \Pic(U)$.

Dealing with such obstruction, I prove a general convergence result
for \'etale cohomology: The canonical map 
$\cH^n_\et(X,\shF)\ra H^n_\et(X,\shF)$ is bijective
for any abelian sheaf $\shF$
provided each $n$-tuple of points $x_1,\ldots,x_n\in X$ 
admits an affine open neighborhood.
This generalizes a result of Artin \cite{Artin 1971},
who assumed that each finite subsets lies in an affine neighborhood.
For noetherian 
schemes such that each pair of points admits an affine open neighborhood,
my result implies that $\bBr(X)=H^2_\et(X,\GG_m)$. 
Here $\bBr(X)$ is Taylor's  
\emph{bigger Brauer group}, defined as 
the group of equivalence classes of central separable
algebras.

Furthermore, we shall see that $H^2_\zar(X,\GG_m)\subset\bBr(X)$
holds for any normal noetherian scheme.
This applies to the nonseparated example  constructed
in \cite{Edidin; Hassett; Kresch; Vistoli 1999},
showing that there are  central separable algebras
neither equivalent to Azumaya algebras nor given by \v{C}ech cocycles.

The paper is organized as follows.
The first section contains observation on
tuples $x_1,\ldots,x_n\in X$ admitting affine open neighborhoods.
In Section 2, I prove the convergence result on \'etale cohomology.
In the next section, I describe the obstruction map
$H^2(X,\shF)\ra\cH^1(X,\shH^1\shF)$ in terms of gerbes and
torsors. The result is purely formal and holds for any site.
Section 4 contains the generalization of Raeburn's and Taylor's
result on the bigger Brauer group. In Section 5, I show that each
Zariski gerbe on a normal noetherian scheme lies in the bigger
Brauer group. The last two sections  contain examples:
Section 6 deals with the nonseparated surface 
from \cite{Edidin; Hassett; Kresch; Vistoli 1999},
and   Section 7 with the proper surfaces without ample
line bundles from \cite{Schroeer 1999}.

\begin{acknowledgements}
I  thank James Borger for   stimulating  discussions,
the Department of Mathematics of the  Massachusetts Institute of
Technology for its hospitality, and
the Deutsche
Forschungsgemeinschaft for financial support.
Moreover, I thank the referee for his comments, which
helped to clarify the paper.
\end{acknowledgements}

\section{Tuples with affine open neighborhoods}

Given a scheme $X$ and an integer $n\geq 2$, we may ask whether each
$n$-tuple $x_1,\ldots,x_n\in X$ 
admits an affine open neighborhood. Such conditions are related
to the existence of ample line bundles
(the generalized Chevalley Conjecture
\cite{Kleiman 1966}, page 327), 
embeddings into toric varieties
\cite{Wlodarczyk 1993}, and \'etale cohomology
\cite{Artin 1971}. In this section, I collect some elementary
results concerning  such conditions.

\begin{proposition}
\mylabel{separated}
Let $X$ be a scheme such that each pair $x_1,x_2\in X$ admits an affine open
neighborhood. Then $X$ is separated.
\end{proposition}

\proof
Let $U_\alpha\subset X$ be the family of all affine open subsets.
Each point in $X\times X$ lies in some subset of the form
$\Spec(\kappa(x_1)\otimes\kappa(x_2))$ with $x_1,x_2\in X$. 
Consequently, the $U_\alpha^2\subset X^2$ form an affine open  covering.
Clearly, the diagonal $\Delta:X\ra X^2$ is a closed embedding over
each $U_\alpha^2$, hence a closed embedding. 
In other words, $X$ is separated.
\qed

\medskip
Given an integer $n\geq 1$ and an $n$-tuple 
$x_1,\ldots,x_n\in X$, consider the subspace 
$S=\Spec(\O_{X,x_1})\cup\ldots\cup\Spec(\O_{X,x_n})$,
which comprises all $x\in X$ specializing to one of the $x_i$.
Setting $\O_S=i^{-1}(\O_X)$, where $i:S\ra X$ is the canonical inclusion,
we obtain a locally ringed space $(S,\O_S)$.
It is covered by the schemes $\Spec(\O_{X,x_i})$.
This covering, however, is not necessarily an open covering,
and  $(S,\O_S)$ is not necessarily a scheme.

\begin{proposition}
\mylabel{scheme}
With the preceding notation, the locally ringed space $(S,\O_S)$
is an affine scheme if the tuple $x_1,\ldots,x_n\in X$ admits an affine open
neighborhood.
\end{proposition}

\proof
To verify this we may assume that $X$ is itself affine.
Now the statement follows form \cite{AC 1-4}, Chap.\ II, \S 3, No.\ 5,
Proposition 17.
\qed

\medskip
I suspect that the converse holds as well.
This is indeed the case under some additional assumptions:

\begin{proposition}
\mylabel{finite type}
Suppose $X$ is  separated and 
of finite type over some noetherian ring $R$.
Then $(S,\O_S)$ is an affine scheme if and only if 
$x_1,\ldots,x_n\in X$
admits an affine open neighborhood.
\end{proposition}

\proof
We already saw that the condition is sufficient and have to verify
necessity. Suppose $(S,\O_S)$ is an affine scheme.
To find the desired affine open neighborhood, we may 
assume that $X$ is reduced by \cite{EGA I}, Corollary 4.5.9.
Adding the generic points $\eta\in X-S$ to the tuple 
$x_1,\ldots,x_n\in X$,
we may also assume that $S\subset X$ is dense.

Choose finitely many sections $g_1,\ldots,g_m\in\Gamma(S,\O_S)$ 
so that the corresponding map
$g:S\ra \AA_R^m$ is injective. 
We may view the $g_i$ as rational functions on $X$
whose domain of definition contains $S$.
Therefore we can replace $X$ by some suitable dense open subset
and assume that the $g_i$ extend to global sections
$f_i\in\Gamma(X,\O_X)$.
In turn, we have a morphism $f:X\ra \AA_R^m$.

Let $U\subset X$ be the subset of $x\in X$ that are isolated in their fiber
$f^{-1}(f(x))$. This is an open subset by 
Chevalley's Semicontinuity Theorem
(\cite{EGA IVc}, Corollary 13.1.4).
By construction, no $x\in S$ admits a generization in $f^{-1}(f(x))$, 
so $S\subset U$.
Replacing $X$ by $U$, we may assume that $f:X\ra \AA^m_R$ has discrete 
fibers. In other words, $f$ is quasifinite.
According to Zariski's Main Theorem 
(\cite{EGA IVc}, Corollary 8.12.6),  there is an open embedding
of $X$ into an affine scheme, hence $\O_X$ is ample.
By \cite{EGA II}, Corollary 4.5.4,
the tuple $x_1,\ldots,x_n\in X$ admits an affine open neighborhood.
\qed

\medskip
Here is another result in this direction.
Recall that a scheme $X$ is called \emph{divisorial}
if the open subset of the form $X_s\subset X$,
where $s$ is a global section of an invertible $\O_X$-module $\shL$,
generate the topology of $X$.
This notion is due to Borelli \cite{Borelli 1963}.

\begin{proposition}
\mylabel{divisorial}
Suppose $X$ is a divisorial noetherian scheme.
Then $(S,\O_S)$ is an affine scheme if and only if $x_1,\ldots,x_n\in X$
admits an affine open neighborhood.
\end{proposition}

\proof
Suppose $(S,\O_S)$ is an affine scheme.
As in the previous proof, we may assume that 
$X$ is reduced and that $S\subset X$ is dense.
By quasicompactness, there is a finitely generated subgroup
$P\subset\Pic(X)$ such that the open subsets $X_s\subset X$, where $s$
ranges over the global sections of the $\shL\in P$, generate the topology.
Choose generators $\shL_1,\ldots,\shL_m\in P$. Then each 
$\shL_i|_S$ is trivial because $S$
is a semilocal affine scheme.
Shrinking $X$ if necessary, we may assume that each $\shL_i$
is trivial. Then $\O_X$ is ample,
and \cite{EGA II}, Corollary 4.5.4
ensures that  $x_1,\ldots,x_n\in X$ admits an affine open neighborhood.
\qed

\section{Obstructions against \v{C}ech cocycles}

Given a scheme $X$, let $X_\et$ be the site of \'etale $X$-schemes.
Its Grothendieck topology is given by the quasicompact \'etale
surjections.
We call such morphism \emph{refinements}, or \emph{\'etale coverings}.
For each abelian sheaf $\shF$
on $X_\et$, we have  cohomology groups
$H^p_\et(X,\shF)$. Sometimes we prefer to deal with
the
\v{C}ech cohomology groups $\cH^p_\et(X,\shF)$ instead.
These groups are related by a natural transformation
$\cH^p_\et(X,\shF)\ra H^p_\et(X,\shF)$ of $\partial$-functors.

For $q\geq 0$, let $\shH^q\shF$ be the presheaf 
$U\mapsto H^q_\et(U,\shF)$. As explained in \cite{Milne 1980},
Chapter III, Proposition 2.7, the composite functor
$\Gamma(X,\shF)=\cH^0(X,\shH^0\shF)$  gives a spectral sequence
$$
\cH^p_\et(X,\shH^q\shF)\Longrightarrow H^{p+q}_\et(X,\shF).
$$
We may view the \v{C}ech cohomology groups
$\cH^p_\et(X,\shH^q\shF)$ with $q>0$ as obstructions against
bijectivity of  
$\cH^p_\et(X,\shF)\ra H^p_\et(X,\shF)$.
The goal of this section is to prove the following vanishing result:

\begin{theorem}
\mylabel{vanishing}
Suppose $X$ is a  noetherian scheme. Let $n\geq 0$ be an integer
such that each $n$-tuple $x_1,\ldots, x_n\in X$  admits an  affine
open neighborhood. Then
$\cH^p_\et(X,\shH^q\shF)=0$ for all
$p<n$, all $q>0$, and any abelian sheaf $\shF$ on $X_\et$.
\end{theorem}

In the case $n=1$, this specializes to the well-known fact 
that $\cH^0_\et(X,\shH^q\shF)=0$ for $q>0$.
The case  $n=\infty$, 
that is, each finite subset lies in an affine open neighborhood,
is Artin's result \cite{Artin 1971}, Corollary 4.2.
We may view  Theorem \ref{vanishing} as a quantitative refinement
of Artin's result. Here is an immediate application:

\begin{corollary}
\mylabel{n-cohomology}
Suppose $X$ is a  noetherian scheme. Let $n\geq 0$ be 
such that each $n$-tuple $x_1,\ldots, x_n\in X$  admits an  affine
open neighborhood. Then the canonical map
$\cH^p_\et(X,\shF)\ra H^p_\et(X,\shF)$ is bijective for $p\leq n$, 
and injective for $p=n+1$. 
\end{corollary}

\proof
The spectral sequence $\cH^p_\et(X,\shH^q\shF)\Rightarrow
H^{p+q}_\et(X,\shF)$ has 
$E_r^{pq}=0$ for all $ p<n$, all $q>0$, and all $r>0$ by Theorem
\ref{vanishing}. Hence the inclusion $E_\infty^{p0}\subset\Gr
H^p_\et(X,\shF)$ is bijective for $p\leq n$. Furthermore
$E_2^{p0}=E_\infty^{p0}$ for $ p\leq n+1$. In turn, the edge map
$\cH^p_\et(X,\shF)\ra H^p_\et(X,\shF)$ is bijective for $p\leq n$,
and injective for $p=n+1$.
\qed

\medskip
Let me also point out the following special case:

\begin{corollary}
\mylabel{toric}
Let $R$ be a  noetherian ring, 
$Y=\bigcup_{\sigma\in\Delta}\Spec(R[\sigma^\vee\cap M])$  a toric
variety, and $X\subset Y$ a subscheme. Then the  map
$\cH^2_\et(X,\shF)\ra H^2_\et(X,\shF)$ is bijective.
\end{corollary}

\proof
According to \cite{Wlodarczyk 1993}, page 709,  each pair of points
in a toric variety admits an affine open neighborhood.
Now the statement follows from Corollary \ref{n-cohomology}.
\qed

\medskip
The proof of Theorem \ref{vanishing} requires a little 
preparation.
Recall that a  scheme is  called \emph{strictly local} if it is the
spectrum of a henselian local ring with separably closed residue
field.

\begin{proposition}
\mylabel{strictly local}
Let $X$ be a quasicompact scheme.
The following are equivalent:
\renewcommand{\labelenumi}{(\roman{enumi})}
\begin{enumerate}
\item  We have $H^p_\et(X,\shF)=0$ for all abelian sheaves
$\shF$ and all $p>0$.
\item Each \'etale covering $U\ra X$ admits a section.
\item The scheme $X$ is affine, and its connected components are
strictly local.
\end{enumerate}
\end{proposition}

\proof
According to \cite{Artin 1971}, Proposition 3.1,  condition (ii)
implies that $X$ is affine. Now the equivalence
(ii)$\Leftrightarrow$(iii)  follows from
\cite{Artin 1971}, Proposition 3.2.
To see the implication (ii)$\Rightarrow$(i), note that each $\shF$-torsor
is trivial on some \'etale covering $U\ra X$, hence trivial, so
the global section functor $H^0(X,\shF)$ is exact.

It remains to verify (i)$\Rightarrow$(ii).
Seeking a contradiction, we assume that some \'etale covering
$f:U\ra X$ admits no  section.
Consider the sheaf $\shF=f_!(\ZZ_U)$. This is the subsheaf 
$f_!(\ZZ_U)\subset f_*(\ZZ_U)$ defined via
extension-by-zero. The \v{C}ech complex for the covering $U\ra X$
is given by
$$
H^0_\et(X,\shF)\stackrel{d_0}{\lra} H^0_\et(U,\shF)
\stackrel{d_1}{\lra} H^0_\et(U^2,\shF).
$$
The constant section $1_U\in H^0_\et(U,f_*\ZZ_U)$ clearly lies in the
subgroup 
$H^0_\et(U,f_!(\ZZ_U))$.
By construction,   $1_U\in H^0_\et(U,\shF)$ lies in the kernel of
$d_1$, but not in the image of $d_0$, and this holds true on all
refinements of $U$. We conclude $\cH^1_\et(X,\shF)\neq 0$.
Since the canonical map $\cH^1_\et(X,\shF)\ra H^1_\et(X,\shF)$ is
injective,  we also have $H^1_\et(X,\shF)\neq 0$, contradiction.
\qed

\medskip
Conforming with \cite{Artin 1971}, Section 3, 
we call a scheme $X$ \emph{acyclic}
if it satisfies the equivalent conditions in Proposition
\ref{strictly local}.
For a point  $x\in X$, let $\O_{X,x}^\sh$ 
be the corresponding \emph{strictly local ring},
that is, the strict henselization of $\O_{X,x}$.
The following is a reformulation of
Artin's fundamental result in \cite{Artin 1971}:

\begin{proposition}
\mylabel{acyclic}
Let $x_1,\ldots, x_n\in X$ be a tuple of points admitting an
affine open neighborhood. Then
the scheme 
$\Spec(\O_{X,x_1}^\sh)\times_X\ldots\times_X\Spec(\O_{X,x_n}^\sh)$ 
is acyclic.
\end{proposition}

\proof
To check this, we may assume that $X$ itself is affine.
Now the assertion follows from \cite{Artin 1971}, Theorem 3.4.
\qed

\medskip
The following improvement will be the key step in proving Theorem
\ref{vanishing}:

\begin{proposition}
\mylabel{refinement}
Suppose $X$ is a noetherian scheme such that
every $(p+1)$-tuple of points in
$X$ admits an affine open neighborhood.
Let 
$U$ be a quasicompact
\'etale $X$-scheme, and $\beta\in H_\et^q(U^{p+1},\shF)$, $q>0$. Let
$V_0,\ldots, V_k$ be quasicompact \'etale
$U$-schemes, and
$x_{k+1}\ldots,x_p\in U$ be points for some $0\leq k\leq p$.
Then there are refinements $V'_i\ra V_i$ for $0\leq i\leq k$, and
affine \'etale neighborhoods $V_i'\ra U$ of $x_i$ for
$k+1\leq i\leq p$, such that
$\beta|_{V'_0\times\ldots\times V'_p}=0$. 
\end{proposition}

\proof
First, we prove by induction on $k$ the following auxiliary
statement:
There are refinements $V'_i\ra V_i$ for $i=0,\ldots,k$
such that 
$
\beta|_{V'_0\times\ldots\times V'_k\times
Z_{k+1}\times\ldots\times Z_p }=0 
$.
Here we write $Z_i=\Spec(\O^\sh_{U,x_i})$ for the
strictly local scheme corresponding to the points $x_i\in U$.

The inductions starts with $k=-1$.
Then  there are no $V_i$, and the assertion boils down 
to  Proposition \ref{acyclic}.
To see this, write each $Z_i=\invlim S_{i,\alpha_i}$ 
as inverse limits
of affine \'etale $X$-schemes $S_{i,\alpha_i}$.
Then $Z_0\times\ldots\times Z_p=\invlim
(S_{0,\alpha_0}\times\ldots\times S_{p,\alpha_p})$, and
\cite{SGA 4b}, Expos\'e VII, Corollary 5.8 tells us that the
canonical map
$$
\dirlim H_\et^q(S_{0,\alpha_0}\times\ldots\times
S_{p,\alpha_p},\shF)\lra 
H_\et^q(Z_0\times\ldots\times Z_p,
\shF_\infty)=0
$$
is bijective, where $\shF_\infty$ is the inverse image of $\shF$.
We conclude that $\beta|_{V'_0\times\ldots\times V'_p}=0$ for
suitable
$V_i'=S_{i,\alpha_i}$. Note that this is the only
step in the   proof where we need the  assumption about affine
neighborhoods of $(p+1)$-tuples.

Now suppose the statement is already true for $k-1$.
Fix  a point $x_k\in V_k$, set
$Z_k=\Spec(\O^\sh_{V_k,x_k})$, and choose refinements
$V'_i\ra V_i$ for
$i=0,\ldots,k-1$ so that 
$\beta|_
{V'_0\times\ldots\times V'_{k-1}\times
Z_k\times\ldots\times Z_p}=0$.
Write $Z_k=\invlim S_\alpha$ as the inverse limit 
of affine \'etale $V_k$-schemes $S_\alpha$.
According to \cite{SGA 4b}, Expos\'e VII, Corollary 5.8, the
canonical map
\begin{gather*}
\dirlim H_\et^q(V'_0\times\ldots\times V'_{k-1}\times
S_\alpha\times Z_{k+1}\times\ldots\times
Z_p,\shF_\alpha)\\ \qquad\qquad\qquad\lra 
H_\et^q(V'_0\times\ldots\times V'_{k-1}\times
Z_k\times Z_{k+1}\times\ldots\times Z_p,
\shF_\infty)
\end{gather*}
is bijective, where $\shF_\alpha$ and $\shF_\infty$  and the inverse
images of $\shF$.
We conclude that $\beta|_{V'_0\times\ldots\times V'_{k-1}\times
S_\alpha\times Z_{k+1}\times\ldots\times
Z_p}=0$ for some suitable index $\alpha$.
If  $S_\alpha\ra V_k$ is surjective,
we are done by setting $V_k'=S_\alpha$. Otherwise, we
finish the argument by applying noetherian
induction to $V_k$. This proves the auxiliary statement.

It remains to construct the desired affine \'etale neighborhoods
$V_i'\ra U$ of the points $x_i\in U$ for $i=k+1,\ldots, p$. 
For this, we write $Z_{k+1}=\invlim T_\alpha$ as the inverse limit 
of affine \'etale $U$-schemes $T_\alpha$. 
Again by \cite{SGA 4b}, Expos\'e VII, Corollary 5.8, the canonical
map
\begin{gather*}
\dirlim H_\et^q(V'_0\times\ldots\times V'_k\times
T_\alpha\times Z_{k+2}\times\ldots\times
Z_p,\shF_\alpha)\\ \qquad\qquad\qquad\lra 
H_\et^q(V'_0\times\ldots\times V'_k\times
Z_{k+1}\times Z_{k+2}\times\ldots\times Z_p,\shF_\infty)
\end{gather*}
is bijective, where $\shF_\alpha$ and $\shF_\infty$  and the inverse
images of $\shF$. As above, we conclude that
$\beta|_{V'_0\times\ldots\times V'_k\times T_\alpha\times
Z_{k+2}\times\ldots\times Z_p}=0$ for some suitable index $\alpha$.
To finish the proof, set $V_{k+1}'=T_\alpha$ and apply induction on
$p-k$.
\qed

\begin{remark}
\mylabel{same refinement}
If there are repetitions among the $V_i$ or the $x_i$,
say $V_i=V_j$ or $x_i=x_j$, then we may also assume   $V_i'=V_j'$, 
by replacing both $V_i'$ and $V_j'$ by $V_i'\times_U V_j'$.
\end{remark}

\medskip
\emph{Proof of Theorem \ref{vanishing}:}
Throughout, we regard $X$ as   base scheme and
products of $X$-schemes as fibered products over $X$.
Fix a \v{C}ech   class $\gamma\in \cH^p(X,\shH^q\shF)$ 
with  $p<n$ and $q>0$.
Choose a refinement $U\ra X$ and a cocycle
$\beta\in H^q(U^{p+1},\shF)$
representing $\gamma$.

It suffices to find a refinement $W\ra U$
with $\beta|_{W^{p+1}}=0$. For this, we shall construct by induction
on $m$   sequences of affine \'etale $U$-schemes
$V_{m,1},\ldots, V_{m,m}$ such that
$
\beta|_{V_{m,i_0}\times\ldots\times V_{m,i_p}}=0
$
for any set of indices $0\leq i_0,\ldots,i_p\leq m$.
This clearly implies $\beta|_{W_m^{p+1}}=0$, where 
$W_m=V_{m,1}\amalg\ldots\amalg V_{m,m}$. In each stage of the
induction,  
$V_{m+1,i}$ will be a refinement of $V_{m,i}$ for
$i=1,\ldots,m$. The induction stops if
$W_m\ra U$ is surjective. We then set $W=W_m$ and have
$\beta|_{W^{p+1}}=0$.

Suppose we already  have constructed
$W_m=V_{m,1}\amalg\ldots\amalg V_{m,m}$ as above, and that
$W_m\ra U$ is not yet surjective.
Fix a point $x\in U$ not in the image and set
$Z=\Spec(\O_{U,x}^\sh)$.  According to
Proposition
\ref{refinement} and Remark \ref{same refinement}, there 
is an affine \'etale neighborhood $V'_{m,m+1}\ra U$
of the point
$x$ such that $\beta|_{{V'}_{m,m+1}^{p+1}}=0$.
Next, fix a tuple of indices
$0\leq i_0,\ldots,i_{p}\leq m+1$.
Applying Proposition \ref{refinement} again, we may replace
the $V_{m,i}'$ for $1\leq i\leq m+1$ by further refinements so that
$
\beta|_{V'_{m,i_0}\times\ldots\times V'_{m,i_{p}}}=0
$.
Since there
are only finitely many such tuples of indices, we may repeat this
inductively until   
$\beta|_{V'_{m,i_0}\times\ldots\times V'_{m,i_{p}}}=0$ 
holds for all  
$0\leq i_0,\ldots,i_{p}\leq m+1$.
Then we set $V_{m+1,i}=V'_{m,i}$ for $i=1,\ldots,m$, and
$V_{m+1,m+1}=V'_{m,m+1}$, and 
$W_{m+1}=V_{m+1,1}\amalg\ldots\amalg V_{m+1,m+1}$.

By construction, the image of $W_{m+1}\ra U$ is strictly larger than
the image of $W_m\ra U$. Using noetherian
induction, we conclude that the mapping $W_m\ra U$ becomes surjective
for some
$m\geq 1$. Hence $W=W_n$ is the desired refinement with 
$\beta|_{W^{p+1}}=0$.
\qed

\section{Gerbes and 2-cohomology}

Theorem \ref{n-cohomology} implies that the  injection   
$\cH^2_\et(X,\shF)\ra H^2_\et(X,\shF) $ is bijective
for any scheme such that each pair $x_1,x_2\in X$ admits an affine open neighborhood.
There is no reason, however, that this   holds in general.
In this section we shall describe the obstruction  in geometric
terms. 

We shall work in an abstract
setting: Fix an arbitrary site with terminal object $X$ and an
abelian sheaf $\shF$. Then we have cohomology
groups $H^p(X,\shF)$.
The spectral sequence
$\cH^p(X,\shH^q\shF)\Rightarrow H^{p+q}(X,\shF)$ gives  an  exact
sequence
$$
0\lra \cH^2(X,\shH^0\shF)\lra H^2(X,\shF)\lra \cH^1(X,\shH^1\shF)
\overset{d}{\lra} \cH^3(X,\shH^0\shF).
$$
The \emph{obstruction map} 
$H^2(X,\shF)\ra \cH^1(X,\shH^1\shF)$ is the obstruction for
a cohomology class to come from  a \v{C}ech cocycle. 
The task now is to describe an obstruction map in terms of gerbes
and torsors.

To do so, let me recall the following geometric  interpretation of
the  universal $\partial$-functor $H^p(X,\shF)$ for $p=0,1,2$:
We may define $H^1(X,\shF)$ as  the group of isomorphism classes of
$\shF$-torsors, and $H^2(X,\shF)$ as the group of
equivalence classes of $\shF$-gerbes. Recall that a \emph{gerbe} is
a stack in groupoids  $\shG\ra X_\et$  satisfying the  following
properties: The objects in $\shG$ are locally isomorphic, and
for each $V\ra X$ there is a refinement $U\ra V$ with $\shG_U$
nonempty.  An \emph{$\shF$-gerbe} is a gerbe $\shG$, together with 
isomorphisms 
$\rho_T:\shF_U\ra \shAut_{T/U}$ for each object
$T\in\shG_U$, such that the $\rho_T$ are compatible with
restrictions, and that the diagram
$$
\begin{CD}
\shF_U @>\rho_T>> \shAut_{T/U}\\
@V\id VV @VV f\mapsto gfg^{-1} V\\
\shF_U @>>\rho_{T'}> \shAut_{T'/U}
\end{CD}
$$
is commutative for each $U$-isomorphism $g:T\ra T'$ (see
\cite{Giraud 1971}, Chapter IV, Definition 2.2.1).
Two $\shF$-gerbes $\shG,\shG'$   are \emph{equivalent} if there
is a functor of stacks $\shG\ra\shG'$ compatible with the
$\shF$-action on automorphism groups. Such functors are
  automatically equivalences  by
\cite{Giraud 1971}, Chapter IV, Corollary 2.2.7.

The $H^p(X,\shF)$, $p=0,1,2$ form a $\partial$-functor as follows:
Given a short exact sequence
$$
0\lra\shF'\lra\shF\lra\shF''\lra 0
$$
and an $\shF''$-torsor $\shT''$, its  liftings 
$(\shT,\shT\ra\shT'')$
to an $\shF$-torsor $\shT$ form
an $\shF'$-gerbe representing the coboundary $\partial(\shT'')$.
According to \cite{Giraud 1971}, Chapter III, Proposition 3.5.1, and
Chapter IV, Lemma 3.4.3, the group
$H^p(X,\shF)$ vanishes on injective sheaves for $p=1,2$,
hence is a universal $\partial$-functor, which justifies the 
notation.

It is easy to express the obstruction map 
$H^2(X,\shF)\ra \cH^1(X,\shH^1\shF)$ in terms of gerbes and torsors:
Let $\shG$ be an $\shF$-gerbe. Choose a covering $U\ra X$ admitting an object
$T\in\shG_U$. Then the sheaf $\Isom(p_0^*T,p_1^*T)$ is an 
$\shF_{U^2}$-torsor on $U^2$, where
$p_i:U^2\ra U$ are the projections omitting the $i$-th factor.
Its isomorphism class is 
a \v{C}ech 1-cochain in $ C^1(U,\shH^1\shF)$.

\begin{lemma}
\mylabel{obstruction torsor}
The $\shH^1\shF$-valued 1-cochain $\Isom(p_0^*T,p_1^*T)$ is a 1-cocycle. 
\end{lemma}

\proof
Set
$\shT=\Isom(p_0^*T,p_1^*T)$, and let 
$p_i:U^3\ra U^2$ be the projections omitting the $i$-th factor.
 We have to see that 
$p_1^*\shT$ is isomorphic to the 
\emph{contracted product} $p_2^*\shT\wedge^\shF p_0^*\shT$.
The latter is the   quotient of $p_2^*\shT\times p_0^*\shT$ 
by the $\shF_{U^2}$-action
$(h_0,h_2)\cdot f=( h_0\circ f,f^{-1}\circ h_2)$.
Using the semisimplicial identities $p_i\circ p_j=p_{j-1}\circ p_i$, $i<j$,
we obtain
\begin{gather*}
p_0^*\shT\simeq \Isom( (p_0p_0)^*T,(p_1p_0)^*T),\quad
p_2^*\shT\simeq \Isom( (p_1p_0)^*T,(p_1p_1)^*T),\\
p_1^*\shT\simeq \Isom( (p_0p_0)^*T,(p_1p_1)^*T).
\end{gather*}
Composition gives a map $p_2^*\shT\times p_0^*\shT\ra p_1^*\shT $, 
which induces
the desired bijection $p_2^*\shT\wedge^\shF p_0^*\shT\simeq p_1^*\shT$.
Note that this bijection is canonical.
\qed

\begin{lemma}
\mylabel{additive}
There is a well-defined linear map 
$H^2(X,\shF)\ra\cH^1(X,\shH^1\shF)$ given by 
$\shG\mapsto\Isom(p_0^*T,p_1^*T)$.
\end{lemma}

\proof
You easily check that the cohomology class of 
$\Isom(p_0^*T,p_1^*T)$ neither depends on the choice
of the refinement $U\ra X$ nor on the choice of the object
$T\in\shG_U$. If $\shG,\shG'$ are two $\shF$-gerbes representing the
same cohomology class, then there is a functor $\shG\ra\shG'$
compatible with the $\shF$-action on automorphism groups.
It follows that the isomorphism class of $\Isom(p_0^*T,p_1^*T)$
depends only on the equivalence class of $\shG$.

It remains to check that the map
$H^2(X,\shF)\ra\cH^1(X,\shH^1\shF)$ is linear. To see this,
choose an injective resolution $\shF\ra\shI^\bullet$.
Given a section $s\in H^0(X,\shI^2)$ contained in the image of
$f:\shI^1\ra\shI^2$, let
$f^{-1}(s)\subset\shI^1$ be the induced $\shI^0/\shF$-torsor, and
$\shG'$ the corresponding
$\shF$-gerbe of $\shI^0$-liftings of $f^{-1}(s)$. 
Let $\shG\subset\shG'$ be the subcategory of liftings $\shI^0_U\ra
f^{-1}(s)_U$ to the trivial
torsor. Since $\shI^0$ is injective, any $\shI^0_U$-torsor is
trivial. Therefore, the inclusion
$\shG\subset\shG'$ is actually a substack hence an equivalence
of $\shF$-gerbes.
Note that any cohomology class is representable by  such an
$\shF$-gerbe
$\shG$, because $\shF\ra\shI^\bullet$ is an injective resolution.

Now choose lifting $\tilde{s}\in H^0(U,\shI^1)$ of $s$ over some
refinement
$U\ra X$. 
This defines the lifting $\shI^0_U\ra f^{-1}(s)_U$,
$0\mapsto \tilde{s}_U$, that is, an object
$T\in\shG_U$.
Now a morphism $p_0^*T\ra p_1^*T$ is precisely a lifting
of $p_1^*(\tilde{s})-p_0^*(\tilde{s})\in H^0(U^2,\shI^0/\shF)$ to
$\shI^0$. Consequently, the torsor
$\Isom(p_0^*T,p_1^*T)$ is nothing but the image of
$p_1^*(\tilde{s})-p_0^*(\tilde{s})\in H^0(U^2,\shI^0/\shF)$ under the
coboundary 
$H^0(U^2,\shI^0/\shF)\ra H^1(U^2,\shF)$ induced by the exact
sequence
$0\ra\shF\ra\shI^0\ra\shI^0/\shF\ra 0$. Using this description, we
immediately infer that
$\shG\mapsto\Isom(p_0^*T,p_1^*T)$ is linear.
\qed

\begin{proposition}
\mylabel{edge map}
An $\shF$-gerbe $\shG$ lies in the image of 
$\cH^2(X,\shF)\ra H^2(X,\shF)$ if and only if the class of 
$\Isom(p_0^*T,p_1^*T)$ vanishes in $\cH^1(X,\shH^1\shF)$. In other
words, we have an exact sequence 
$0\ra \cH^2(X,\shF)\ra H^2(X,\shF)\ra\cH^1(X,\shH^1\shF)$. 
\end{proposition}

\proof
According to \cite{Giraud 1971}, Chapter IV, Corollary 2.5.3,
an $\shF$-gerbe $\shG$ comes from $\cH^2(X,\shF)$ if and only
if
it admits an object
$T\in\shG_U$ over some refinement $U\ra X$ with 
$p_0^*(T)\simeq p_1^*(T)$, hence $\Isom(p_0^*T,p_1^*T)$ is trivial.

Now suppose $\shT=\Isom(p_0^*T,p_1^*T)$ has trivial
cohomology class. Replacing $U$ by a refinement, we find an
$\shF$-torsor   
$\shP$ on $U$ with 
$\Isom(p_1^*\shP,p_0^*\shP)\simeq\shT$.
According to \cite{Giraud 1971}, Chapter III, Proposition 2.3.2
there is a twisted object $T'\in\shG_U$ satisfying $\shP=\Isom(T,T')$.
Then $\Isom(p_0^*T',p_1^*T')$, being isomorphic to
$$
\Isom(p_0^*T',p_0^*T)\wedge
\Isom(p_0^*T,p_1^*T)\wedge
\Isom(p_1^*T,p_1^*T')=p_0^*(\shP^{-1})\wedge\shT \wedge
p_1^*(\shP),
$$
is trivial, and we conclude that the class of $\shG$ lies in 
$\cH^2(X,\shF)$.
\qed

\section{Central separable algebras}
\mylabel{central separable algebra}

In this section I apply Theorem \ref{vanishing}
to   the bigger Brauer group. Throughout, $X$ denotes a noetherian
scheme. Let me recall some notions from Raeburn and Taylor 
\cite{Raeburn; Taylor 1985}.  Given two coherent
$\O_X$-modules $\shE,\shF$ and a pairing $\lambda:\shF\otimes\shE\ra\O_X$, 
we obtain a coherent
$\O_X$-algebra $\shE\otimes^\lambda\shF$ as follows: 
The underlying $\O_X$-module is $\shE\otimes\shF$,
and the multiplication law is
$$
(e\otimes f)\cdot (e'\otimes f') = e\lambda(f,e')\otimes f' =
e\otimes \lambda(f,e')f'.
$$
Usually, $\shE\otimes^\lambda\shF$ is neither commutative nor unital. 
We are mainly interested in the case that $\lambda$ is surjective;
this ensures that $\shE$, $\shF$, and 
$\shE\otimes^\lambda\shF$ are faithful $\O_X$-modules.

Now let $\shA$ be a coherent $\O_X$-algebra. A \emph{splitting} for $\shA$
is a quadruple $(\shE,\shF,\lambda,s)$, where $\shE,\shF$ are coherent 
$\O_X$-modules,
$\lambda:\shF\otimes\shE\ra\O_X$ is a surjective pairing, and 
$s:\shA\ra\shE\otimes^\lambda\shF $ is an $\O_X$-algebra bijection.
We say that $\shA$ is \emph{elementary} if it admits a splitting.
If there is an \'etale covering $U\ra X$ so that $\shA_U$ admits a splitting,
we say that $\shA$ is a \emph{central separable} algebra.

Suppose $\shA$ is a central separable  algebra. For each \'etale map $U\ra X$,
let $\shS_U$ be the groupoid of splittings for $\shA_U$; a morphism
$(\shE,\shF,\lambda,s)\ra(\shE',\shF',\lambda',s')$ of splittings is 
a pair of bijections $e:\shE\ra\shE'$ and $f:\shF\ra\shF'$ such that
the diagrams
$$
\begin{CD}
\shF\otimes\shE @>\lambda >> \O_X\\
@Vf\otimes e VV @VV\id V\\
\shF'\otimes\shE' @>\lambda' >> \O_X
\end{CD}
\quadand\quad
\begin{CD}
\shA @>s>> \shE\otimes\shF\\
@V\id VV @VVe\otimes fV\\
\shA @>s'>> \shE'\otimes\shF'\\
\end{CD}
$$
commute.
Clearly, the fibered category $\shS\ra X_\et$ is a stack
in Giraud's sense (\cite{Giraud 1971}, Chapter II, Definition
1.2.1).  According
to \cite{Raeburn; Taylor 1985}, Lemma 2.3, the splittings for $\shA$
are locally isomorphic.  Furthermore, each splitting 
$(\shE,\shF,\lambda,s)$
comes along with  a sheaf homomorphism
$$
\GG_m\lra\shAut_{(\shE,\shF,\lambda,s)},\quad 
\xi\longmapsto(\xi,1/\xi),
$$
which is bijective by \cite{Raeburn; Taylor 1985}, Lemma 2.4.
In other words, $\shS$ is a $\GG_m$-gerbe.
So  each central separable  algebra 
$\shA$ defines via the gerbe $\shS$ a
cohomology class in $ H^2_\et(X,\GG_m)$.

Next, let us recall Taylor's definition of the  bigger Brauer group.
You easily check that central separable  algebras are closed
under taking opposite algebras and tensor products.
Two central separable  algebras $\shA,\shA'$ are  called
\emph{equivalent} if there are elementary algebras $\shB,\shB'$ with 
$\shA\otimes\shB\simeq\shA'\otimes\shB'$.
The set of equivalence classes $\bBr(X)$ is called the 
\emph{bigger Brauer group}.
Addition is given by tensor product, and inverses are given
by opposite algebras.

The map $\shA\mapsto\shS$ induces an inclusion 
$\bBr(X)\subset H^2_\et(X,\GG_m)$
of abelian groups. 
Raeburn and Taylor \cite{Raeburn; Taylor 1985} showed that
this inclusion is a bijection provided that  each finite subset
of $X$ admits a common affine neighborhood.
We may relax this assumptions:

\begin{theorem}
\mylabel{bigger brauer}
Let $X$ be a noetherian scheme with the property that
each pair $x,y\in X$ admits an affine open neighborhood.
Then $\bBr(X)=H^2_\et(X,\GG_m)$.
\end{theorem}

\proof
The proof of Raeburn and Taylor actually shows that, on
an arbitrary noetherian scheme,  each \v{C}ech 2-cohomology class
comes from a coherent central separable $\O_X$-algebra
(\cite{Raeburn; Taylor 1985}, Theorem 3.6).
According to Theorem \ref{vanishing}, we
have $\cH^2_\et(X,\GG_m)=H^2_\et(X,\GG_m)$, 
and in turn $\bBr(X)=H^2_\et(X,\GG_m)$.
\qed

\section{Normal noetherian schemes}

Hilbert's Theorem 90 implies that the  map
$H^2_\zar(X,\O_X^\times)\ra H^2_\et(X,\GG_m)$ 
is injective. The goal of this section is to  construct central
separable algebras representing classes from this subgroup.
Throughout, we shall assume that $X$ is a normal noetherian scheme.

Let $\shDiv_X$ and $\shZ^1_X$ be the sheaves   
of Cartier divisors and
Weil divisors with respect to the Zariski topology, and
$\shP_X=\shZ^1_X/\shDiv_X$ the corresponding quotient sheaf.
Similarly, let $\Div(X)$ and
$Z^1(X)$ be the groups of Cartier divisors and Weil divisors, and
$\Cl(X)=Z^1(X)/\Div(X)$. Setting $P(X)=\Gamma(X,\shP_X)$,  we obtain
an inclusion
$\Cl(X)\subset P(X)$.

\begin{proposition}
\mylabel{zariski}
Let $X$ be a normal noetherian scheme. 
Then there is a canonical identification
$H^2_\zar(X,\O_X^\times)=P(X)/\Cl(X)$.
\end{proposition}

\proof
Let $\shM_X^\times$ be the sheaf of invertible rational functions.
The exact sequence
$
1\ra \O_X^\times\ra \shM_X^\times\ra\shDiv_X\ra 0
$
gives an exact sequence
$$
H^1_\zar(X,\shM_X^\times)\lra H^1_\zar(X,\shDiv_X)
\overset{\partial}{\lra} H^2_\zar(X,\O_X^\times)\lra
H^2_\zar(X,\shM_X^\times).
$$
The outer groups $H^n_\zar(X,\shM_X^\times)$ vanish;
to check this, use the spectral
sequence $H^p_\zar(X,R^qi_*\O^\times_{X^{(0)}})\Rightarrow
H^{p+q}_\zar(X^{(0)},\O^\times_{X^{(0)}})$,  where $i: X^{(0)}\ra X$
is the inclusion of the generic points. Now the exact sequence
$
0\ra \shDiv_X\ra\shZ^1_X\ra \shP_X\ra 0
$
gives an exact sequence
$$
\Div(X)\lra Z^1(X)\lra P(X)\overset{\partial}{\lra}
H^1_\zar(X,\shDiv_X)\lra H^1_\zar(X,\shZ^1_X).
$$
The term on the right vanishes, because  $\shZ^1_X$ is flabby, and
the result follows.
\qed

\medskip
Weil divisors give rise to central separable algebras in 
the following way:
Given finitely many $C_1,\ldots,C_n\in
Z^1(X)$, consider the coherent reflexive sheaves
$$
\shE=\bigoplus_{\nu=1}^n\O_X(C_\nu)
\quadand
\shF=\bigoplus_{\nu=1}^n\O_X(-C_\nu).
$$
Let $\lambda_{\nu\mu}:\O_X(C_\nu)\otimes\O_X(-C_\mu)\ra\O_X$ be the
pairing defined as
$$
f\otimes g\longmapsto 
\left\{
\begin{array}{l}
 f(g)\quad \text{if $\nu=\mu$,}\\
 0 \quad \text{otherwise.}
\end{array}
\right.
$$
The $(n\times n)$-matrix of pairings $\lambda=(\lambda_{\nu\mu})$
defines a pairing 
$\lambda:\shF\otimes\shE\ra\O_X$.
As described in Section \ref{central separable algebra}, this yields
a coherent $\O_X$-algebra
$\shA=\shE\otimes^\lambda\shF$.

Clearly, the pairing $\lambda:\shF\otimes\shE\ra\O_X$ is surjective
if at each point
$x\in X$  at least one Weil divisor $C_i$ is Cartier.
Under this assumption, $\shA$ is a central separable
$\O_X$-algebra endowed with a splitting.
We shall use such algebras for the following result:

\begin{theorem}
\mylabel{inclusion}
Suppose  $X$ is a normal noetherian scheme. Then we have inclusions 
$H^2_\zar(X,\O_X^\times)\subset \bBr(X)$ of subgroups
in $H^2_\et(X,\GG_m)$.
\end{theorem}

\proof
Fix a class $\alpha\in H^2_\zar(X,\O_X^\times)$, 
and choose a representant $s\in P(X)$ with respect
to the canonical surjection $P(X)\ra H^2_\zar(X,\O_X^\times)$ from
Proposition \ref{zariski}. 
Then $s\in P(X)$ is given by a collection of Weil
divisors
$D_i\in Z^1(U_i)$ on some open covering $X=U_1\cup\ldots\cup U_n$, 
such that $D_i-D_j$ are Cartier
on the overlaps $U_{ij}=U_i\cap U_j$.
We may extend each $D_i$ from $U_i$ to $X$ and denote the resulting
Weil divisor
$D_i\in Z^1(X)$ by the same letter.
For each $U_i\subset X$, set
$$
\shE_i=\bigoplus_{\nu=1}^n\O_{U_i}(D_i-D_\nu)
\quadand
\shF_i=\bigoplus_{\nu=1}^n\O_{U_i}(D_\nu-D_i).
$$
As above, this yields a coherent $\O_{U_i}$-algebra
$\shA_i=\shE_i\otimes^{\lambda_i}\shF_i$. They are central
separable because
$D_i-D_\nu$ is Cartier on $U_i$ for $\nu=i$.

These   $\O_{U_i}$-algebras glue together as follows:
For each overlap $U_{ij}=U_i\cap U_j$, consider the invertible
$\O_{U_{ij}}$-module
$\shL_{ij}=\O_{U_{ij}}(D_i-D_j)$. We have canonical isomorphisms
$$
\shE_i|_{U_{ij}}\otimes\shL_{ji}\lra\shE_j|_{U_{ij}}
\quadand
\shL_{ij}\otimes\shF_i|_{U_{ij}}\lra\shF_j|_{U_{ij}}.
$$
The canonical bijections $\shL_{ji}\otimes\shL_{ij}\ra\O_{U_{ij}}$ 
yield
isomorphisms $\lambda_{ij}:\shA_j|_{U_{ij}}\ra\shA_i|_{U_{ij}}$.
These isomorphisms obviously satisfy the cocycle condition
$\lambda_{ij}\circ\lambda_{jk}=\lambda_{ik}$ on
triple overlaps.
We deduce that there is a coherent central separable $\O_X$-algebra
$\shA$ with $\shA|_{U_i}=\shA_i$.

It remains to check that the $\O_X^\times$-gerbe $\shS$ of splittings
for
$\shA$ has cohomology class $\alpha\in H^2_\zar(X,\O_X^\times)$.
Let $f:\shZ^1_X\ra\shP_X$ be the canonical surjection, and   $\shG'$
be the $\O_X^\times$-gerbe of $\shM_X^\times$-liftings for the
$\shDiv_X$-torsor
$f^{-1}(s)\subset\shZ^1_X$. Then $\shG'$ has class $\alpha\in
H^2_\zar(X,\O_X^\times)$ because $s\mapsto\alpha$.
Note that $H^1(U,\shM_X^\times)=0$ for any open subset $U\subset X$.
Therefore, 
the fibered 
subcategory $\shG\subset\shG'$ of liftings of $f^{-1}(s)$
to the trivial $\shM_X^\times$-torsor $\shM_X^\times$
is an $\O_X^\times$-subgerbe.

To finish the proof, we construct a functor $\shG\ra\shS$
compatible with  $\O_X^\times$-actions.
Suppose we have an object in $\shG$ over an open subset
$V\subset X$, that is, an
equivariant map $\shM_V^\times\ra f^{-1}(s)|_V$.
Let $D\in\Gamma(V,f^{-1}(s))$ be the image of the unit section
$1\in\Gamma(V,\shM_X^\times)$. Then $D-D_i$ are Cartier on 
$V_i=V\cap U_i$.
Consider the coherent reflexive $\O_{V_i}$-modules
$\shE_i'=\shE_i\otimes\O_{V_i}(D-D_i)$ and
$\shF_i'=\shF_i\otimes\O_{V_i}(D_i-D)$. We have splittings
$\shA|_{V_i}=\shE'_i\otimes^{\lambda'_i}\shF'_i$.
Note that
$$
\shE'_i=\bigoplus_{\nu=1}^n\O_{V_i}(D-D_\nu)
\quadand
\shF'_i=\bigoplus_{\nu=1}^n\O_{V_i}(D_\nu-D).
$$
Obviously, the sheaves $\shE'_i$ glue together and give
a coherent $\O_V$-module $\shE'$. Similarly, the $\shF'_i$ glue
and give a
coherent
$\O_V$-module $\shF$. In turn, we obtain a splitting
$\shA|_V=\shE'\otimes^{\lambda'}\shF'$.

Summing up, we have defined for each   object in
$\shG$ an object in $\shS$. It is easy to see that
this construction is functorial
and respects the $\O_X^\times$-action on automorphism groups.
Therefore, the central separable $\O_X$-algebra $\shA$
has class $\alpha\in H^2_\zar(X,\O_X^\times)$.
\qed

\medskip
Next, we describe the obstruction against cocycles.
Fix a cohomology class $\alpha\in H^2_\zar(X,\O_X^\times)$ and choose
$s\in P(X)$ mapping to
$\alpha$. Then there is an open covering $X=U_1\cup\ldots\cup U_n$
and Weil divisors $D_i\in Z^1(U_i)$ representing $s|_{U_i}$, such
that $D_i-D_j$ are Cartier on the overlaps $U_{ij}$.

\begin{proposition}
\mylabel{zariski obstruction}
The cocycle
$U_{ij}\mapsto\O_{U_{ij}}(D_i-D_j)$ represents the image of $\alpha$
under the obstruction map
$H^2_\zar(X,\O_X^\times)\ra \cH^1_\zar(X,\shH^1(\O_X^\times))$.
\end{proposition}

\proof
Consider the exact sequence
$
1\ra\O_X^\times\ra\shM_X^\times\ra \shZ^1\ra \shP_X\ra 0.
$
Since $H^1_\zar(U,\shM_X^\times)=0$ for any open subset $U\subset X$,
we may argue as in the proof of Proposition \ref{additive} and infer
that $U_{ij}\mapsto\O_{U_{ij}}(D_i-D_j)$ represents
the image of $\alpha$ under the obstruction map 
$H^2_\zar(X,\O_X^\times)\ra\cH^1_\zar(X,\shH^1(\O_X^\times))$.
\qed

\section{Nonseparated surfaces}

Recall that the Brauer group $\Br(X)\subset H^2_\et(X,\GG_m)$ is the
subgroup generated by Azumaya algebras, and that  the
\emph{cohomological Brauer group} 
$\Br'(X)\subset H^2_\et(X,\GG_m)$ is the torsion subgroup.
In this section we discuss  the example of
Edidin, Hassett, Kresch, and Vistoli \cite{Edidin; Hassett;
Kresch; Vistoli 1999} of a scheme with $\Br(X)\neq\Br'(X)$.  

Let $A$
be a strictly local normal noetherian ring of dimension two that is
nonfactorial. In other words, $A$ is neither regular nor an
$E_8$-singularity (\cite{Brieskorn 1968}, Proposition 3.3). 
Set $Y=\Spec(A)$
and let
$W\subset Y$ be the complement of the closed point.
Define $X=U_1\cup U_2$ as the union of two copies of $Y$ 
glued along $W$.
Then $X$ is a normal nonseparated surface with two closed
points $x_1\in U_1,x_2\in U_2$.

The theory
of quotient stacks was used in 
\cite{Edidin; Hassett; Kresch; Vistoli
1999} to prove $\Br(X)\neq \Br'(X)$.  Let us present a different
argument. The covering $X=U_1\cup U_2$ gives an exact sequence
$$
\bigoplus_{i=1}^2 H^1(U_i,\GG_m)
\lra H^1(U_1\cap U_2, \GG_m) \lra H^2(X,\GG_m) \lra
\bigoplus_{i=1}^2 H^2(U_i,\GG_m)
$$
for both Zariski and \'etale cohomology.
The outer terms vanish because the  $U_i$ are strictly local.
Together with Hilbert's Theorem  90, 
this implies $H^2_\et(X,\GG_m)= H^2_\zar(X,\GG_m)$.
Hence \emph{every cohomology class comes from a central separable
$\O_X$-algebra} by Theorem \ref{inclusion}.

Using Proposition \ref{zariski} and the canonical 
bijections $\Cl(X)=\Pic(W)=\Cl(Y)$, we conclude
$$
H^2_\et(X,\GG_m)=\Cl(U_1)\oplus\Cl(U_2)/\Cl(Y) \simeq \Cl(Y) \neq 0.
$$
This 
implies $\cH^2_\zar(X,\O_X^\times)=0$. 
Indeed, suppose some class in $H^2_\zar(X,\GG_m)$ represented by
a pair of Weil divisors
$(D_1,D_2)\in\Cl(U_1)\oplus\Cl(U_2)$ vanishes in the obstruction
group
$\cH^1_\zar(X,\shH^1(\O_X^\times))$. By Proposition
\ref{zariski obstruction}, the invertible sheaf $\O_W(D_1-D_2)$ is of
the form
$\shL_1|_W\otimes\shL_2|_W$ with $\shL_i\in\Pic(U_i)=0$. 
It follows that our pair $(D_1,D_2)$ is zero in $H^2_\et(X,\GG_m)$.
Summing up,
\emph{only the trivial cohomology class comes from a cocycle}.

As explained in \cite{Schroeer 2001},
Proposition 1.5,
each Azumaya $\O_X$-algebra $\shA$ is of the form $\shEnd(\shE)$ for some
reflexive $\O_X$-module $\shE$, say of rank $r>0$,  with
$\shE_{U_i} = \O_{U_i}(D_i)^{\oplus r}$ for some Weil divisors
$D_i\in Z^1(U_i)$. Furthermore, the class of $\shA$
is the image of $-(D_1,D_2)$ in $H^2_\zar(X,\GG_m)$.
Since $\Gamma(U_1,\shA)=\Gamma(W,\shA)=\Gamma(U_2,\shA)$, 
we have $D_1\sim D_2$ and 
conclude  $\Br(X)=0$. In other words, \emph{only the trivial
cohomology class comes from an Azumaya algebra}.

\section{Nonprojective proper surfaces}

In this section I discuss the cohomology groups
$H^2_\et(X,\GG_m)$ for some nonprojective proper surfaces constructed
in \cite{Schroeer 1999}. Let me recall the construction:
Fix an algebraically closed ground field
$k$, let $E$ be an elliptic curve, 
and choose two closed points $e_1,e_2\in E$.
Let $Y\ra \PP^1\times E$ be the blowing-up of 
the points $(0, e_1),(\infty,e_2)$, 
and $g:Y\ra X$ the contraction of
the strict transforms $E_1,E_2\subset Y$ of $0\times E, \infty\times
E$. Then $X$ is a proper normal algebraic surface containing
two  singularities $x_1,x_2\in X$ of genus $g$. 
As explained in \cite{Schroeer 1999}, it has no ample line bundles if
the divisor classes $e_1,e_2\in\Pic(E)\otimes\QQ$ are linearly independent.

\begin{proposition}
\mylabel{group}
We have $H^2_\zar(X,\O_X^\times)\simeq\Pic(E)/\ZZ e_1+\ZZ e_2$.
\end{proposition}

\proof
The sheaf $\shP_X=\shZ^1_X/\shDiv_X$ is a skyscraper sheaf supported
by the singular locus $\left\{x_1,x_2\right\}$, with stalks
$\shP_{x_i}=\Cl(\O_{X,x_i})$. 
According to Proposition \ref{zariski},
we have
$$
H^2_\zar(X,\O_X^\times)=\Cl(\O_{X,x_1})\oplus\Cl(\O_{X,x_2})/\Cl(X).
$$
The terms on the right are
$\Cl(\O_{X,x_i})=\Pic(Y\otimes\O_{X,x_i})/\ZZ E_i$, where
$Y\otimes\O_{X,x_i}$ denotes the fiber product
$Y\times_X\Spec(\O_{X,x_i})$. Moreover, the canonical mapping
$\Pic(Y\otimes\O_{X,x_i})\ra \Pic(Y\otimes\O^\wedge_{X,x_i})$ is
injective. Grothendieck's Existence Theorem gives
$\Pic(Y\otimes\O^\wedge_{X,x_i})=\Pic(nE_i)$ for some
$n>0$, and we have $\Pic(nE_i)=\Pic(E_i)$ because $E_i$ is elliptic.
Consequently
$\Cl(\O_{X,x_i})=\Pic(E_i)/\ZZ e_i$

The group  $\Cl(X)$ is generated by the images of
$\Pic(E)$, $\Pic(\PP^1)$, and the exceptional divisors
for the contraction
$Y\ra\PP^1\times E$. The latter two types restrict to zero in
$\Cl(\O_{X,x_i})$. The result now follows from the snake lemma.
\qed

\begin{proposition}
\mylabel{bijective}
The  inclusion  $H^2_\zar(X,\O_X^\times)\subset H^2_\et(X,\GG_m)$ is
bijective.
\end{proposition}

\proof
We have $H^2_\et(E,\GG_m)=0$ because the ground field is algebraically
closed (\cite{GB}, Corollary 1.2). In turn
$H^2_\et(\PP^1\times E,\GG_m)$ vanishes
(\cite{Gabber 1981}, page 193, Theorem 2).
By birational invariance,  $H^2_\et(Y,\GG_m)$ vanishes as well
(\cite{GB}, Corollary 7.2).
Now the commutative diagram
$$
\begin{CD}
\Pic(Y_\zar)@>>>H^0_\zar(X,R^1g_*\O_Y^\times)@>>>
H^2_\zar(X,\O_X^\times)@>>>H^2_\zar(Y,\O_Y^\times)\\
@VVV @VVV @VVV @VVV\\
\Pic(Y_\et)@>>>H^0_\et(X,R^1g_*\GG_m)@>>>H^2_\et(X,\GG_m)@>>>
H^2_\et(Y,\GG_m).
\end{CD}
$$
The map on the left is bijective by Hilbert's Theorem 90.
The map next to the left  is nothing but the sum of the maps
$\Pic(Y\otimes\O_{X,x_i})\ra \Pic(Y\otimes\O^\sh_{X,x_i})$.
But both $\Pic(Y\otimes\O_{X,x_i})$ and 
$\Pic(Y\otimes\O^\sh_{X,x_i})$
are  equal to $\Pic(E_i)$ as shown  in the proof for Proposition
\ref{group}. We infer
$H^2_\zar(X,\GG_m)=H^2_\et(X,\GG_m)$ using the 5-Lemma.
\qed

\begin{proposition}
\mylabel{no cocycle}
We have $\cH^2_\zar(X,\O_X^\times)=0$.
\end{proposition}

\proof
We have to check 
that the   map 
$H^2_\zar(X,\GG_m)\ra \cH^1_\zar(X,\shH^1\GG_m)$ is injective.
Pick some $s\in P(X)$.
Choose an open covering $U_i\subset X$ so that 
$s$ lifts to Weil divisors  $D_i\in Z^1(U_i)$.
The image of $s$ in $\cH^1_\zar(X,\shH^1\GG_m)$ is represented by the
1-cocycle
$U_{ij}\mapsto \O_{U_{ij}}(D_i-D_j)$.
Suppose this class is zero. After refining the covering,
there are Cartier divisors $C_i\in\Div(U_i)$ with $D_i-D_j=C_i-C_j$.
After re-indexing, we may assume $x_1\in U_1$ and $x_2\in U_2$.
Since $D_1$ is principal on $\Spec(\O_{X,x_1})$, and $D_2$ is principal 
on $\Spec(\O_{X,x_2})$, we infer that
$C_1-C_2$ is a principal divisor on the Dedekind scheme
$S=\Spec(\O_{X,x_1})\times_X \Spec(\O_{X,x_2})$, 
which comprises all points $x\in X$ 
with $\left\{x_1,x_2\right\}\subset\overline{\left\{x\right\}}$. 
But this implies that $s$ is the restriction of
a global reflexive rank one sheaf, such that 
$s$ maps to zero  in $H^2_\zar(X,\GG_m)$.
\qed

\begin{question}
Is the inclusion $\cH^2_\et(X,\GG_m)\subset H^2_\et(X,\GG_m)$
bijective?  Does the obstruction group $\cH^1_\et(X,\shH^1\GG_m)$
vanish?
\end{question}


\end{document}